\documentclass{amsart}
\usepackage{amsmath,amssymb,amsbsy,amsfonts,latexsym,
amsthm,amscd,amsxtra,amsrefs,cite}
\usepackage[T1]{fontenc}
\usepackage{graphics,latexsym,enumerate}

\newtheorem{lemma}{Lemma}
\newtheorem{thm}{Theorem}

\newtheorem{cor}{Corollary}

\begin{document}

\title{A sharp Subelliptic Sobolev Embedding Theorem with weights}
\author{Po-Lam YUNG}
  \address{Department of Mathematics \\ the Chinese University of Hong Kong \\ Ma Liu Shui \\ Shatin \\ Hong Kong}
  \email{plyung@math.cuhk.edu.hk}
  \thanks{The author was supported by a Titchmarsh Fellowship at the University of Oxford, a Junior Research Fellowship at St. Hilda's College, and a direct grant for research 4053120 at the Chinese University of Hong Kong.}

\begin{abstract}
The purpose of this short article is to prove some potential estimates that naturally arise in the study of subelliptic Sobolev inequalites for functions. This will allow us to prove a local subelliptic Sobolev inequality with the optimal amount of smoothing, as well as a variant of that which describes quantitatively an improvement of the inequality as one gets away from certain characteristic varieties.
\end{abstract}

\maketitle

\section{Statement of results}

Subelliptic Sobolev-type estimates in general have received a lot of attention over the years. We list some results that share a similar theme as ours: Capogna-Danielli-Garofalo \cite{MR1312686}, Cohn-Lu-Wang \cite{MR2345338}, Franchi-Gallot-Wheeden \cite{MR1314734}, Franchi-Lu-Wheeden \cite{MR1343563}, \cite{MR1354890}, \cite{MR1383947}, Franchi-P{\'e}rez-Wheeden \cite{MR1744780}, Franchi-Wheeden \cite{MR1455453}, Jerison \cite{MR850547}, Lu \cite{MR1202416}, \cite{MR1286482}, Lu-Wheeden \cite{MR1642822}, \cite{MR1792599}, Muckenhoupt - Wheeden \cite{MR0340523}, P{\'e}rez-Wheeden \cite{MR1818113}, \cite{MR1962949} and Sawyer-Wheeden \cite{MR1175693}.

In \cite{PL}, the author has proved a Sobolev inequality for the $\overline{\partial}_b$-complex on $(0,q)$ forms on a certain class of CR manifolds of finite type. In this current work, the focus will be on functions (rather than forms), and the result is real-variable in nature.

To describe our results, we need to introduce some notations. Following Nagel, Stein and Wainger \cite{MR793239} and \cite{MR1882665}, let $\Omega \subseteq \mathbb{R}^N$ be a connected open set, and let $Y_1, \dots, Y_q$ be a list, possibly with repetitions, of smooth real vector fields on $\Omega$. Assume that to each $Y_j$ we associate an integer $d_j \geq 1$, called the formal degree of $Y_j$. The collection $\{Y_j\}_{j=1}^q$ is said to be of finite homogeneous type on $\Omega$ if they span $\mathbb{R}^N$ at every point in $\Omega$, and that for each $1 \leq j,k \leq q$, $$[Y_j,Y_k] = \sum_{d_l \leq d_j + d_k} c_{j,k}^l (x) Y_l$$ for some $c_{j,k}^l \in C^{\infty}(\Omega)$. For example, if $X_1, \dots, X_n$ are smooth real vector fields on $\Omega$ that satisfy Hormander's condition, meaning that successive commutators of $X_1,\dots,X_n$ of some finite length $r$ already span the tangent space at every point of $\Omega$, then if $\{Y_j\}$ is the collection of successive commutators of $X_1,\dots,X_n$ up to length $r$, it is of finite homogeneous type. 

With such a collection $\{Y_j\}$, one can then define a control metric $\rho$ as follows. For each $\delta > 0$, let $C(\delta)$ be the set of absolutely continuous curves $\phi \colon [0,1] \to \Omega$ such that $\phi'(t) = \sum_{j=1}^q a_j(t) Y_j(\phi(t))$ with $|a_j(t)| \leq \delta^{d_j}$ for all $j$ and almost all $t \in [0,1]$. For $x,y \in \Omega$, let $\rho(x,y) = \inf \{ \delta > 0 \colon \text{there is a curve } \phi \in C(\delta) \text{ such that } \phi(0)=x \text{ and } \phi(1)=y \}.$ We shall write $B(x,\delta)$ for the metric ball centered at $x$ and of radius $\delta$, namely $\{y \in \Omega \colon \rho(x,y) < \delta\}$, and $V(x,y)$ for the Lebesgue measure of the ball $B(x,\rho(x,y))$.

If now $I$ is an $N$-tuple $(i_1,\dots,i_N)$, $1 \leq i_j \leq q$, we write $$Q_I = \sum_{j=1}^N d_{i_j},$$ and $$\lambda_I(x) = \det(Y_{i_1}, \dots, Y_{i_N})(x)$$ for $x \in \Omega$. Here we are taking the determinant of the $N \times N$ matrix, whose $j$-th column is the component of $Y_{i_j}$ in the coordinate basis $\frac{\partial}{\partial x_1}, \dots, \frac{\partial}{\partial x_N}$. These numbers are important in computing the volumes of the metric balls (see Theorem~\ref{thm:HT} below). It is also in terms of these numbers that we state our main results.

\begin{thm}\label{thm:WPE}
For each $N$-tuple $I$ and each compact subset $E$ of $\Omega$, the map $$T_If(x) = |\lambda_I(x)|^{\frac{1}{Q_I}}  \int_E \frac{\rho(x,y)}{V(x,y)} f(y) dy$$ maps $L^p(E)$ boundedly into $L^{p^*}(E)$, where $$\frac{1}{p^*} = \frac{1}{p} - \frac{1}{Q_I}, \quad 1 < p < Q_I.$$ It also maps $L^1(E)$ into weak-$L^{\frac{Q_I}{Q_I-1}}(E)$. Here $dy$ is the Lebesgue measure on $E$, and all the $L^p$ spaces are taken with respect to the Lebesgue measure on $E$.
\end{thm}

We also have:

\begin{thm}\label{thm:mPE}
For each $N$-tuple $I$ and each compact subset $E$ of $\Omega$, the map $$Tf(x) = \int_E \frac{\rho(x,y)}{V(x,y)} f(y) dy$$ maps $L^1(E,dy)$ boundedly into weak-$L^{\frac{Q_I}{Q_I-1}}(E,d\mu_I)$, where $$d\mu_I(x) := |\lambda_I(x)|^{\frac{1}{Q_I-1}} dx,$$ and $dx$ is the Lebesgue measure on $E$.
\end{thm}

These allow us to prove the following subelliptic Sobolev inequality for Hormander's vector fields:

\begin{thm}\label{thm:WSI}
Let $X_1, \dots, X_n$ be smooth real vector fields on a connected open set $\Omega \subseteq \mathbb{R}^N$, whose commutators of length $\leq r$ span at every point of $\Omega$. List the commutators of $X_1, \dots, X_n$ of length $\leq r$ as $Y_1, \dots, Y_q$, and define $\lambda_I(x)$ for each $N$-tuple $I$ and $x \in \Omega$ as above. Let $\Omega'$ be a relatively compact open subset of $\Omega$ with smooth boundary and $I$ be an $N$-tuple. Then for each $f \in C^{\infty}(\overline{\Omega'})$, we have 
$$
\left(\int_{\Omega'} |f(x)|^{p^*} |\lambda_I(x)|^{\frac{p}{Q_I-p}} dx\right)^{\frac{1}{p^*}} \leq C \left( \int_{\Omega'}  |\nabla_b f(x)|^p + |f(x)|^p dx \right)^{\frac{1}{p}},
$$
where 
\begin{equation} \label{eq:pp*}
\frac{1}{p^*} = \frac{1}{p}-\frac{1}{Q_I} \quad \text{and} \quad 1 \leq p < Q_I.
\end{equation}
\end{thm}

Here the length of the subelliptic gradient $|\nabla_b f|$ is defined by
$$
|\nabla_b f|^2 := |X_1 f|^2 + \dots + |X_n f|^2.
$$
By picking $I$ to be the $N$-tuple with minimal $Q_I$ such that $|\lambda_I| \simeq 1$ around each point in $\overline{\Omega'}$, and patching the estimates together, we obtain the following corollary:

\begin{cor} \label{cor:SIQ}
Let $X_1, \dots, X_n$ be as in Theorem~\ref{thm:WSI}. For each $x \in \Omega$, let $Q(x)$ be the non-isotropic dimension at $x$, defined by 
$$
Q(x) := \sum_{j=1}^r j n_j(x), \qquad n_j(x) := \dim V_j(x) - \dim V_{j-1}(x).
$$
where $V_j(x)$ is the span of the commutators of $X_1, \dots, X_n$ of length $\leq j$ at $x$. Let $\Omega'$ be a  relatively compact open subset of $\Omega$ with smooth boundary, and define the non-isotropic dimension $Q$ of $\overline{\Omega'}$ by setting $$Q := \sup_{x \in \overline{\Omega'}} Q(x).$$ Then for any $f \in C^{\infty}(\overline{\Omega'})$, 
\begin{equation}\label{eq:SSI}
\|f\|_{L^{p^*}(\Omega')} \leq C (\|\nabla_b f\|_{L^p(\Omega')} + \|f\|_{L^p(\Omega')}),
\end{equation}
where $$\frac{1}{p^*} = \frac{1}{p} - \frac{1}{Q}, \quad 1 \leq p < Q.$$
\end{cor}

This is a subelliptic Sobolev inequality with a maximal degree of smoothing. It implies Proposition 1 of \cite{PL}, which was stated there without proof. See also the work of Caponga, Danielli and Garofalo \cite{MR1266765}, Varopoulos \cite{MR1070036} and Gromov \cite[Section 2.3.D'']{MR1421823}. We shall also prove that the exponent $p^*$ given in the corollary is \emph{always} sharp; in other words, this inequality cannot hold for any bigger values of $p^*$. See Section~\ref{sect:cor} below.

The theorem for general $I$, on the other hand, says that one gets more smoothing as soon as one looks at regions where $\lambda_I(x)$ does not degenerate to 0, and it tells us how such an improved inequality degenerates as $\lambda_I(x)$ degenerates to 0.

For instance, if we have the vector fields $\frac{\partial}{\partial x_1}$ and $x_1^{r-1}\frac{\partial}{\partial x_2}$ on $\mathbb{R}^2$, then Theorem~\ref{thm:WSI} (and a rescaling argument) implies that 
$$
\left(\int_{\mathbb{R}^2} |f(x_1,x_2)|^{p^*} |x_1|^{(r-1)\frac{p}{2-p}} dx_1dx_2\right)^{\frac{1}{p^*}} \leq C \left( \left\|\frac{\partial f}{\partial x_1} \right\|_{L^p(\mathbb{R}^2)} + \left\|x_1^{r-1} \frac{\partial f}{\partial x_2} \right\|_{L^p(\mathbb{R}^2)}\right),
$$
where
$$\frac{1}{p^*} = \frac{1}{p} - \frac{1}{2}, \quad 1 \leq p < 2.$$ Note that the factor $|x_1|^{^{(r-1)\frac{p}{2-p}}}$ on the left hand side tends to 0 as the vector field $x_1^{r-1} \frac{\partial}{\partial x_2}$ degenerates, when one moves towards the axis where $\{x_1 = 0\}$.

In the case where the underlying space is a homogeneous group, however, Theorem~\ref{thm:WSI} does not improve upon the known results, because $\lambda_I \equiv 0$ unless $Q_I$ is bigger than or equal to the homogeneous dimension of the group.

The observation (as depicted in Theorem~\ref{thm:WSI} above) that one can use different weights at different points in a Sobolev or isoperimetric inequality is certainly not new; see for example Franchi-Gallot-Wheeden \cite{MR1314734}, Franchi-Lu-Wheeden \cite{MR1343563}, \cite{MR1354890}, \cite{MR1383947}. Typically when one uses different weights, one attach a `dimension' to every point, that may vary not just with the weights being used, but also from point to point. In Franchi-Wheeden \cite{MR1455453}, they introduced the concept of a \emph{compensation couple}, in an attempt to `stablize' the dimensions used at various points (when such a couple exists). 
In light of this, it would be natural to ask what the `best' weight is in any given situation. Unfortunately our results have little to say in this direction.

One can also prove the following variant of Theorem~\ref{thm:WSI}, where instead of a zeroth order term in $f$ on the right hand side, we have $f$ minus the average of $f$ on the left-hand side. Moreover, one can replace the smoothness condition on $\Omega'$, by a weaker Boman chain condition: an open set $\Omega' \subset \Omega$ will be said to satisfy the Boman chain condition $\mathcal{F}(\tau,M)$ for some $\tau \geq 1$, $M \geq 1$, if there exists a covering $W$ of $\Omega'$ by (Carnot-Caratheodory) balls $B$, such that 
$$
\sum_{B \in W} \chi_{\tau B}(x) \leq M \chi_{\Omega'}(x) \quad \text{for all $x \in \Omega$},
$$
and there exists a ``central'' ball $B_1 \in W$, which can be connected to every ball $B \in W$ by a finite chain of balls $B_1, \dots, B_{\ell(B)} = B$ of $W$ so that $B \subset M B_j$ for $j=1,\dots,\ell(B)$, with the additional property that $B_j \cap B_{j-1}$ contains a ball $R_j$ such that $B_j \cup B_{j-1} \subset M R_j$ for $j = 2, \dots, \ell(B)$. (Here all balls are Carnot-Caratheodory balls. Also, $\tau B$ denotes a ball that has the same center as $B$, but $\tau$ times the radius, and $\chi_{S}$ denotes the characteristic function of a set $S$.) 

\begin{thm} \label{thm:WSI2}
Let $\Omega'$ be a relatively compact open subset of $\Omega$ that satisfies the Boman chain condition $\mathcal{F}(\tau,M)$ for some $\tau \geq 1$, $M \geq 1$. For any $N$-tuple $I$ and any $1 \leq p < Q_I$, let $w_{I,p}(x)$ be the weight defined by
$$
w_{I,p}(x) := |\lambda_I(x)|^{\frac{p}{Q_I-p}}.
$$
Assume that $w_{I,p}(x) dx$ is a doubling measure. Then for any Lipschitz functions $f$ on $\overline{\Omega'}$, we have
\begin{equation} \label{eq:fwPIfinal}
\left( \int_{\Omega'} |f(x) - f_{\Omega'}|^{p^*} w_{I,p}(x) dx \right)^{\frac{1}{p^*}} 
\leq  C \left( \int_{\Omega'} |\nabla_b f(x)|^p dx \right)^{\frac{1}{p}},
\end{equation}
where 
\begin{equation} \label{eq:fwav}
f_{\Omega'} := \frac{\int_{\Omega'} f(x) w_{I,p}(x) dx}{\int_{\Omega'} w_{I,p}(x) dx},
\end{equation}
and $p^*$ is as in (\ref{eq:pp*}).
\end{thm}

Note that if $\Omega'$ is a relatively compact subset of $\Omega$ with smooth boundary, then it satisfies a Boman chain condition for some $\tau \geq 1$, $M \geq 1$. More generally, the same is true for all John domains \cite{MR1427074}, \cite{MR1343563}, so the above theorem applies for such $\Omega'$'s as well.

On the other hand, we only managed to establish such a theorem under the additional doubling condition on our weighted measure $w_{I,p}(x) dx$. Such a doubling condition is satisfied by a number of important examples (e.g. the Grushin type example given by the vector fields $\frac{\partial}{\partial x_1}$ and $x^{r-1} \frac{\partial}{\partial x_2}$ on $\mathbb{R}^2$), but could fail when say $\lambda_I(x)$ vanishes on some open set (e.g. if $Y_1 = \frac{\partial}{\partial x_1}$, $Y_2 = \frac{\partial}{\partial x_2}$, $Y_3 = (1-a(x)) \frac{\partial}{\partial x_1} + a(x) \frac{\partial}{\partial x_2}$ on $\mathbb{R}^2$, where $a(x)$ vanishes on some non-trivial open set, then when $I = \{1,3\}$, $\lambda_I(x) = a(x)$ vanishes on some non-trivial open set). It is not clear whether such doubling conditions are really necessary.

It is an interesting question whether the pair of weights $(w_{I,p}(x),1)$ satisfies the local balance condition in the work of Chanillo-Wheeden \cite{MR805809} (i.e. condition (1.5) of \cite{MR1343563}). If it is, then Theorem~\ref{thm:WSI2} would follow from the work of Franchi-Lu-Wheeden in \cite{MR1343563}.

The author thanks the referee for suggesting the possibility of Theorem~\ref{thm:WSI2}, and for raising the above question about the pair of weights $(w_{I,p}(x),1)$. The author is also grateful to the referee for numerous very helpful comments.

\section{Proof of Theorem~\ref{thm:WPE}}

Let $Y_1,\dots,Y_q$ be of finite homogeneous type in $\Omega$ as in the previous section. We recall the following Theorem of Nagel, Stein and Wainger, from \cite{MR793239} and \cite{MR1882665}:

\begin{thm}[Nagel-Stein-Wainger]\label{thm:HT}
Let $E$ be a compact subset of $\Omega$. Then for all $x \in E$ and all $\delta < \text{diam}_{\rho}(E)$, where $\text{diam}_{\rho}(E)$ is the diameter of $E$ with respect to the metric $\rho$, we have $$|B(x,\delta)| \simeq \max_J |\lambda_J(x)| \delta^{Q_J},$$ where the maximum is over all $N$-tuples $J$. (Hereafter we write $\simeq$ or $\lesssim$ when the implicit constants depend only on $E$.) 
\end{thm}

In particular, the Lebesgue measure is doubling on $E$ with respect to the metric balls defined by $\rho$, and $V(x,y) \simeq V(y,x)$ for all $x,y \in E$. 

Now to prove Theorem~\ref{thm:WPE}, fix any $N$-tuple $I$ and a compact subset $E$ of $\Omega$. We observe the following pointwise estimate for the kernel of $T_I$: 
$$|\lambda_I(x)|^{\frac{1}{Q_I}} \frac{\rho(x,y)}{V(x,y)} = \frac{(|\lambda_I(x)| \rho(x,y)^{Q_I})^{\frac{1}{Q_I}}}{V(x,y)} \lesssim V(x,y)^{\frac{1}{Q_I}-1} \simeq V(y,x)^{\frac{1}{Q_I}-1}.$$ This is just a simple consequence of Theorem~\ref{thm:HT}. Hence for any $x \in E$ and any $\alpha > 0$, the set $$\left\{y \in E \colon |\lambda_I(x)|^{\frac{1}{Q_I}} \frac{\rho(x,y)}{V(x,y)} > \alpha \right\} \subseteq  \left\{y \in E \colon V(x,y) \lesssim \alpha^{-\frac{Q_I}{Q_I-1}} \right\},$$ the latter of which is a metric ball centered at $x$, whose Lebesgue measure is $\simeq \alpha^{-\frac{Q_I}{Q_I-1}}$ uniformly in $x$. Similarly, for any $y \in E$, $$\left\{x \in E \colon |\lambda_I(x)|^{\frac{1}{Q_I}} \frac{\rho(x,y)}{V(x,y)} > \alpha \right\} \subseteq  \left\{x \in E \colon V(y,x) \lesssim \alpha^{-\frac{Q_I}{Q_I-1}} \right\},$$ which is a metric ball centered at $y$, and has Lebesgue measure $\simeq \alpha^{-\frac{Q_I}{Q_I-1}}$ uniformly in $y$. Hence $T_I$ maps $L^p(E)$ to weak-$L^{p^*}(E)$ whenever $\frac{1}{p^*} = \frac{1}{p} - \frac{1}{Q_I}$, where $1 \leq p < Q_I$. We now invoke the following version of the Marcinkiewicz interpolation theorem (which can be found, e.g. in Lemma 15.3 of Folland-Stein~\cite{MR0367477}):

\begin{lemma}
Let $k$ be a measurable function on $E \times E$ such that for some $r > 1$, $k(x,\cdot)$ is weak-$L^r$ uniformly in $x$, and $k(\cdot,y)$ is weak-$L^r$ uniformly in $y$. Then the operator $f(x) \mapsto \int k(x,y) f(y) dy$ is bounded from $L^p(E)$ to $L^q(E)$ whenever $$\frac{1}{q} + 1 = \frac{1}{p} + \frac{1}{r}, \quad 1 < p < q < \infty.$$
\end{lemma}

From the above estimates for the kernel of $T_I$, if we apply the lemma with $r = \frac{Q_I}{Q_I-1}$, we see that $T_I$ mapping $L^p(E)$ to $L^{p^*}(E)$, whenever $1 < p < Q_I$. 

\section{Proof of Theorem~\ref{thm:mPE}}

We now turn to the proof of Theorem~\ref{thm:mPE}. Let $r := \frac{Q_I}{Q_I-1}$. Then $r > 1$, and weak-$L^r$ is a normed space. So by Minkowski inequality, if $f \in L^1(E,dy)$, then 
$$
\|Tf\|_{\text{weak-}L^r(E,d\mu_I)} \leq \|f\|_{L^1(E,dy)} \sup_{y \in E} \left\| \frac{\rho(x,y)}{V(x,y)} \right\|_{\text{weak-}L^r(E,d\mu_I(x))}.
$$
Since $d\mu_I(x) = |\lambda(x)|^{\frac{1}{Q_I-1}} dx$, it suffices to show that for any $y \in E$ and $\alpha > 0$,
$$
\int_{\{x \in E \colon \frac{\rho(x,y)}{V(x,y)} > \alpha\}} |\lambda(x)|^{\frac{1}{Q_I-1}} dx \lesssim \alpha^{-r} \quad \text{uniformly in $y$.}
$$
Now $\{x \in E \colon \frac{\rho(x,y)}{V(x,y)} > \alpha\} \subseteq \{x \in E \colon \frac{\rho(y,x)}{V(y,x)} \gtrsim \alpha\}$, and the latter is a metric ball centered at $y$. Let $\delta_{\alpha}$ be its radius, so that it is equal to $B(y,\delta_{\alpha})$; then 
\begin{equation}\label{eq:deltaalpha}
\frac{\delta_{\alpha}}{|B(y,\delta_{\alpha})|} \simeq \alpha.
\end{equation}
Recall that by Theorem~\ref{thm:HT}, $|\lambda_I(x)| \delta_{\alpha}^{Q_I} \lesssim |B(x,\delta_{\alpha})|$. Hence for any $x \in B(y,\delta_{\alpha})$, we have
$$
|\lambda_I(x)| \lesssim |B(x,\delta_{\alpha})| \delta_{\alpha}^{-Q_I} \lesssim |B(y,\delta_{\alpha})| \delta_{\alpha}^{-Q_I}.
$$
(The last inequality follows from the doubling property of the Lebesgue measure with respect to the metric balls.) Hence 
\begin{align*}
\int_{\{x \in E \colon \frac{\rho(x,y)}{V(x,y)} > \alpha\}} |\lambda(x)|^{\frac{1}{Q_I-1}} dx 
&\lesssim \int_{B(y,\delta_{\alpha})} |B(y,\delta_{\alpha})|^{\frac{1}{Q_I-1}} \delta_{\alpha}^{-\frac{Q_I}{Q_I-1}} dx \\
&= |B(y,\delta_{\alpha})|^{\frac{Q_I}{Q_I-1}} \delta_{\alpha}^{-\frac{Q_I}{Q_I-1}} \simeq \alpha^{-\frac{Q_I}{Q_I-1}} = \alpha^{-r},
\end{align*}
the second-to-last equality following from (\ref{eq:deltaalpha}). This completes our proof.

\section{Proof of Theorem~\ref{thm:WSI}}

We can now prove Theorem~\ref{thm:WSI}. First recall the following pointwise potential estimate, versions of which are well-known: Let $\Omega'$ be a relatively compact open subset of $\Omega$ with smooth boundary, and $E = \overline{\Omega'}$ be its closure. For any $f \in C^{\infty}(E)$ and any $x \in E$, we have 
$$
|f(x)| \lesssim \int_{E} \frac{\rho(x,y)}{V(x,y)} (|\nabla_b f(y)| + |f(y)|) dy.
$$
In fact, this estimates follows from an analysis of the fundamental solution of the sum of squares operator $-\sum_{j=1}^n X_j^* X_j$, as was analyzed in Nagel-Stein-Wainger~\cite{MR793239} and S\'anchez-Salle~\cite{MR762360}. (See also discussion following formula (1.2) of Franchi-Lu-Wheeden~\cite{MR1343563}.) Theorem~\ref{thm:WSI} then follows readily from Theorem~\ref{thm:WPE}, in the case when $p > 1$. If $p = 1$, we need a well-known truncation argument, to show that the strong type bound follows from the weak type bound we proved in Theorem~\ref{thm:mPE}; c.f. Long-Rui \cite{MR1187073}, and the exposition in Haj{\l}asz \cite{MR1886617} or Chapter 3 of Heinonen \cite{MR1800917}. The crucial reason why this truncation argument works is that we are not letting the potential operator $T$ act on arbitrary functions; instead they are all acting on some gradient of a single function. We also need the fact that the gradients here are taken using real (rather than complex) vector fields. 

First, according to Theorem~\ref{thm:mPE}, for all $N$-tuples $I$ and all $\alpha > 0$,
$$
\int_{E \cap \{|f|> \alpha\}} |\lambda_I(x)|^{\frac{1}{Q_I-1}} dx \lesssim \alpha^{-\frac{Q_I}{Q_I-1}} (\|\nabla_b f\|_{L^1(E)} + \|f\|_{L^1(E)})^{\frac{Q_I}{Q_I-1}}.
$$
This originally holds for all $f \in C^{\infty}(E)$, but this also holds for all $f \in W^{1,1}(E)$, because one can approximate such functions both in $W^{1,1}(E)$ and almost everywhere by smooth functions in $E$.

Now let $f \in C^{\infty}(E)$, and for any integer $j$ let  
$$
f_j = 
\begin{cases} 
0 &\quad \text{if $|f| \leq 2^{j-1}$} \\
|f|-2^{j-1} &\quad \text{if $2^{j-1} \leq |f| \leq 2^j$}\\
2^{j-1} &\quad \text{if $|f| \geq 2^j$}
\end{cases}.
$$ 
Then $f_j \in W^{1,1}(E)$ (this is a qualitative statement; we will not need bounds on the $W^{1,1}(E)$ norms of the $f_j$'s), 
and
$$
\nabla_b f_j = 
\begin{cases} 
\nabla_b |f| &\quad \text{if $2^{j-1} < |f| < 2^j$}\\
0 &\quad \text{otherwise}
\end{cases}
$$ 
in distribution. 
Hence by the weak-type $(L^1(dy),L^{\frac{Q_I}{Q_I-1}}(d\mu_I))$ result above, 
\begin{align*}
\int_{E \cap \{f_j > \alpha\}} |\lambda_I(x)|^{\frac{1}{Q_I-1}} dx 
&\lesssim \alpha^{-\frac{Q_I}{Q_I-1}} \left( \int_{2^{j-1} < |f| < 2^j} |\nabla_b f(x)| dx + \int_E |f_j(x)|dx \right)^{\frac{Q_I}{Q_I-1}} 
\end{align*}
because wherever $f \ne 0$, $$|\nabla_b |f|| \leq |\nabla_b f|$$ (here we need $X_1, \dots, X_n$ to be real vector fields, because $f$ may be complex-valued). It then follows that 
\begin{align*}
&\int_E |f(x)|^{\frac{Q_I}{Q_I-1}} |\lambda_I(x)|^{\frac{1}{Q_I-1}} dx\\ 
\leq& \sum_{j=-\infty}^{\infty} (2^{j+1})^{\frac{Q_I}{Q_I-1}} \int_{E \cap \{ 2^j < |f| \leq 2^{j+1}\}} |\lambda_I(x)|^{\frac{1}{Q_I-1}} dx  \\
\leq& \sum_{j=-\infty}^{\infty} (2^{j+1})^{\frac{Q_I}{Q_I-1}} \int_{E \cap \{ f_j > 2^{j-1}\}} |\lambda_I(x)|^{\frac{1}{Q_I-1}} dx \\
\lesssim &  \sum_{j=-\infty}^{\infty} (2^{j+1})^{\frac{Q_I}{Q_I-1}} (2^{j-1})^{-\frac{Q_I}{Q_I-1}} \left( \int_{2^{j-1} < |f| < 2^j} |\nabla_b f(x)| dx + \int_E |f_j(x)| dx \right)^{\frac{Q_I}{Q_I-1}} \\
\lesssim & \left(\|\nabla_b f\|_{L^1(E)} + \|f\|_{L^1(E)}\right)^{\frac{Q_I}{Q_I-1}}
\end{align*}
as desired.

\section{Proof of Corollary~\ref{cor:SIQ} and its sharpness} \label{sect:cor}

Let $\Omega'$ be as in the previous section, and $E = \overline{\Omega'}$. Recall that at every point $x \in E$, we defined a local non-isotropic dimension $Q(x)$, and from its definition, it is clear that there exists a neighborhood $U_x$ of $x$ and an $N$-tuple $I_x$ such that the degree of $I_x$ is $Q(x)$, and such that $|\lambda_{I_x}| \simeq 1$ on $U_x \cap E$. From Theorem~\ref{thm:WSI}, it follows that for all $f \in C^{\infty}(U_x \cap E)$ and all $1 \leq p < Q(x)$, we have $$\|f\|_{L^{\frac{Q(x)p}{Q(x)-p}}(U_x \cap E)} \lesssim \|\nabla_b f\|_{L^p(U_x \cap E)} + \|f\|_{L^p(U_x \cap E)}.$$ Since $Q = \sup_{x \in E} Q(x)$, by taking a partition of unity and gluing the estimates, we see that Corollary~\ref{cor:SIQ} follows.

We now prove that the exponent $p^*$ in Corollary~\ref{cor:SIQ} is always the best possible. This will follow from a consideration of an approximate dilation invariance. For this we need to introduce a suitable coordinate system and a non-isotropic dilation near a point $x_0 \in E$ where $Q(x_0) = \sup_{x \in E} Q(x)$. 
Let $x_0$ be as such, and let $\{X_{jk} \colon 1 \leq j \leq r, 1 \leq k \leq n_j\}$ be a collection of vector fields that satisfies the following:
\begin{enumerate}[(a)]
 \item Each $X_{jk}$ is a commutator of $X_1, \dots, X_n$ of length $j$;
 \item For each $1 \leq j_0 \leq r$, $\{X_{jk} \colon 1 \leq j \leq j_0, 1 \leq k \leq n_j\}$ restricts at $x_0$ to a basis of $V_{j_0}(x_0)$.
\end{enumerate}
In particular $$\sum_{j=1}^r jn_j = Q(x_0) = Q.$$
Then for some small $\varepsilon > 0$,
\begin{align}
[-\varepsilon,\varepsilon]^N &\to \mathbb{R}^N \notag \\
u &\mapsto \exp(u \cdot X') x_0 \label{eq:coord}
\end{align} 
defines a normal coordinate system in a neighborhood $U_0$ of $x_0$ in $\mathbb{R}^N$; here $\exp(X)x_0$ is the time-1-flow along the integral curve of the vector field $X$ beginning at $x_0$, and
$$
u \cdot X' = \sum_{j = 1}^r \sum_{k=1}^{n_j} u_{jk} X_{jk}
$$ 
where $u = (u_{jk})_{1 \leq j \leq r, 1 \leq k \leq n_j}$. For simplicity we shall consistently write $u$ for $\exp(u \cdot X')x_0 \in U_0$. This coordinate system allows us to define the associated non-isotropic dilation: for $u = (u_{jk}) \in U_0$ and $\lambda > 0$, write $$\lambda \cdot u := (\lambda^{j} u_{jk})_{1 \leq j \leq r, 1 \leq k \leq n_j}$$ as long as the latter is in $U_0$ (and we leave this undefined if it is not in $U_0$).

Now if $\alpha = (j_1k_1, \dots, j_sk_s)$ is a multiindex, we shall let $u^{\alpha}$ be the monomial $u_{j_1k_1}u_{j_2k_2} \dots u_{j_sk_s}.$ It is said to have non-isotropic degree $|\alpha|=j_1+\dots+j_s$ because $$(\lambda \cdot u)^{\alpha} = \lambda^{|\alpha|} u^{\alpha}.$$ A function $f$ of $u$ is said to vanish to non-isotropic order $l$ at $0$ if its Taylor series expansion consists of terms whose non-isotropic degrees are all $\geq l$. Note that if $$X_l = \sum_{j=1}^{r} \sum_{k=1}^{n_j} a^{l}_{jk}(u) \frac{\partial}{\partial u_{jk}}$$ on $U_0$ for $1 \leq l \leq n$, then by the Campbell-Hausdorff formula, each $a^l_{jk}(u)$ vanishes to non-isotropic order $j-1$ at $u=0$. (c.f. Section 10 of \cite{MR0436223}). In what follows we Taylor expand $a^l_{jk}$ at $0$ such that $$a^l_{jk}(u) = p^l_{jk}(u) + e^l_{jk}(u)$$ where $p^l_{jk}(u)$ are homogeneous polynomials of non-isotropic degree $j-1$ and $e^l_{jk}(u)$ vanish to non-isotropic order $j$ at $0$. Define $$W_l = \sum_{j=1}^{r} \sum_{k=1}^{n_j} p^{l}_{jk}(u) \frac{\partial}{\partial u_{jk}}, \quad E_l = \sum_{j=1}^{r} \sum_{k=1}^{n_j} e^{l}_{jk}(u) \frac{\partial}{\partial u_{jk}}$$ on $U_0$, for $1 \leq l \leq n$.

Given $p \geq 1$, let $q$ be an exponent for which (\ref{eq:SSI}) holds for all $f \in C^{\infty}(E)$. We shall show $q \leq p^*$. In fact then the inequality holds for all $f \in C^{\infty}_c(U_0 \cap E)$ (just extend $f$ by zero to all of $E$): 
$$
\left( \int_{U_0 \cap E} |f(x)|^q dx \right)^{\frac{1}{q}} \lesssim \left( \int_{U_0 \cap E}  |\nabla f(x)|^p + |f(x)|^p  dx \right)^{\frac{1}{p}}.
$$
But we can also parameterize $U_0 \cap E$ by the $u$ coordinates we introduced above, and use the Lebesgue measure $du$ with respect to this $u$ coordinates in place of $dx$ in the above inequality. This is because $du$ is a smooth density times $dx$, and vice versa. Hence for all $f \in C^{\infty}_c(U_0 \cap E)$,
\begin{equation}\label{eq:NLS2}
\left( \int_{U_0 \cap E} |f(u)|^q du \right)^{\frac{1}{q}} \lesssim \left( \int_{U_0 \cap E} |\nabla f(u)|^p + |f(u)|^p  du \right)^{\frac{1}{p}}.
\end{equation}

Now pick an open set $U_1 \subseteq U_0 \cap E$ such that $\lambda \cdot u \in U_0 \cap E$ for all $\lambda \in [0,1]$. This is possible because $E$ is the closure of an open set with smooth boundary. Then take $f \in C^{\infty}_c(U_1)$ that is not identically zero. For each $\delta \in (0,1)$, let 
$$
f_{\delta}(u) := 
\begin{cases}
f(\delta^{-1} \cdot u) &\quad \text{if $\delta^{-1} \cdot u \in U_1$}\\ 
0 &\quad \text{otherwise}
\end{cases}.
$$ 
Applying (\ref{eq:NLS2}) to $f_{\delta}$ in place of $f$, we get 
\begin{align*}
\delta^{\frac{Q}{q}} \|f\|_{L^{q}(U_1)} 
&\leq C \left(\sum_{j=1}^n \|W_j(f_{\delta})\|_{L^p(U_1)} + \sum_{j=1}^n \|E_j(f_{\delta})\|_{L^p(U_1)} + \|f_{\delta}\|_{L^p(U_1)} \right) \\
&= C \sum_{j=1}^n \delta^{-1+\frac{Q}{p}} \|W_j f\|_{L^p(U_1)} +  O(\delta^{\frac{Q}{p}})
\end{align*}
by the homogeneity of the vector fields $W_j$ and $E_j$. Letting $\delta \to 0$, we get $\frac{Q}{q} \geq -1 + \frac{Q}{p}$, i.e. $$\frac{1}{q} \geq \frac{1}{p}-\frac{1}{Q}.$$ Hence $q \leq p^*$ as desired.

We remark that a similar argument shows that Theorem~1 of \cite{PL} cannot hold for any value of $Q$ smaller than the one stated there.

\section{Proof of Theorem~\ref{thm:WSI2}}

To prove Theorem~\ref{thm:WSI2}, an important starting point is a representation formula, as derived in Lu-Wheeden \cite{MR1642822}. It was proved, in Theorem 1 there, that if $B$ is any Carnot-Caratheodory ball in $\Omega'$, then 
\begin{equation} \label{eq:oprepform}
|f(x)-L(f,B)| \leq C \int_{B} \frac{\rho(x,y)}{V(x,y)}  |\nabla_b f(y)| dy, \quad \text{for all $x \in B$},
\end{equation}
where $L(f,B) := |B|^{-1} \int_B f(y) dy$ is the Lebesgue average of $f$ over $B$, and $C$ is an appropriate constant. If $1 < p < Q_I$, it then follows from Theorem~\ref{thm:WPE} that
\begin{equation} \label{eq:wPIBall}
\left( \int_{B} |f(x) - L(f,B)|^{p^*} w_{I,p}(x) dx \right)^{\frac{1}{p^*}} \leq C \left( \int_{B} |\nabla_b f(x)|^p dx \right)^{\frac{1}{p}}.
\end{equation}
The corresponding weak-type $(1,1^*)$ bound, and the truncation argument used in the proof of Theorem~\ref{thm:WSI} shows that the above inequality remains true when $p = 1$. Now we patch these estimates together,
using the Boman chain condition satisfied by $\Omega'$. In fact, since we assumed that $w_{I,p}(x) dx$ is a doubling measure, using Theorem~3.7 of Franchi-Lu-Wheeden \cite{MR1343563}, one then concludes the existence of some constant $A(f,\Omega')$, such that 
\begin{equation} \label{eq:fwPIpre}
\left( \int_{\Omega'} |f(x) - A(f,\Omega')|^{p^*} w_{I,p}(x) dx \right)^{\frac{1}{p^*}} \leq C \left( \int_{\Omega'} |\nabla_b f(x)|^p dx \right)^{\frac{1}{p}}.
\end{equation}
Then it is a standard argument that one can replace $A(f,\Omega')$ by the average $f_{\Omega'}$, where $f_{\Omega'}$ is defined as in (\ref{eq:fwav}). In fact, 
$$
f_{\Omega'} - A(f,\Omega') = \frac{\int_{\Omega'} [f(x) - A(f,\Omega')] w_{I,p}(x) dx}{\int_{\Omega'} w_{I,p}(x) dx}.
$$
Hence by Jensen's inequality,
$$
|f_{\Omega'} - A(f,\Omega')| \leq \left( \frac{\int_{\Omega'} |f(x) - A(f,\Omega')|^{p^*} w_{I,p}(x) dx}{\int_{\Omega'} w_{I,p}(x) dx} \right)^{\frac{1}{p^*}}.
$$
Moving the denominator to the left hand side, and applying (\ref{eq:fwPIpre}), we then see that
\begin{equation} \label{eq:fwPIpre2}
\left( \int_{\Omega'} \left| f_{\Omega'} - A(f,\Omega') \right|^{p^*} w_{I,p}(x) dx \right)^{\frac{1}{p^*}} \leq C \left( \int_{\Omega'} |\nabla_b f(x)|^p dx \right)^{\frac{1}{p}}.
\end{equation}
Combining (\ref{eq:fwPIpre}) and (\ref{eq:fwPIpre2}) finally gives the desired estimate (\ref{eq:fwPIfinal}).

We remark that in place of (\ref{eq:oprepform}), we could also have used formula (1.2) of Franchi-Lu-Wheeden \cite{MR1343563} instead, which states the same inequality as in (\ref{eq:oprepform}), except that the integral on the right hand side was over a ball $cB$ for some $c > 1$ instead of just over $B$. In that case, we would have to replace, in the right hand side of (\ref{eq:wPIBall}), $B$ by $cB$, but Theorem 3.7 of Franchi-Lu-Wheeden \cite{MR1343563} still applies, and we get the same inequality (\ref{eq:fwPIpre}) as desired.

\bibliography{div_curl_paper}{}

\def\cprime{$'$}
\begin{bibdiv}
\begin{biblist}

\bib{MR1427074}{incollection}{
      author={Buckley, S.},
      author={Koskela, P.},
      author={Lu, G.},
       title={Boman equals {J}ohn},
        date={1996},
   booktitle={X{VI}th {R}olf {N}evanlinna {C}olloquium ({J}oensuu, 1995)},
   publisher={de Gruyter, Berlin},
       pages={91\ndash 99},
      review={\MR{1427074 (98m:43013)}},
}

\bib{MR1312686}{article}{
      author={Capogna, Luca},
      author={Danielli, Donatella},
      author={Garofalo, Nicola},
       title={The geometric {S}obolev embedding for vector fields and the
  isoperimetric inequality},
        date={1994},
        ISSN={1019-8385},
     journal={Comm. Anal. Geom.},
      volume={2},
      number={2},
       pages={203\ndash 215},
      review={\MR{1312686 (96d:46032)}},
}

\bib{MR1266765}{article}{
      author={Capogna, Luca},
      author={Danielli, Donatella},
      author={Garofalo, Nicola},
       title={An isoperimetric inequality and the geometric {S}obolev embedding
  for vector fields},
        date={1994},
        ISSN={1073-2780},
     journal={Math. Res. Lett.},
      volume={1},
      number={2},
       pages={263\ndash 268},
      review={\MR{MR1266765 (95a:46048)}},
}

\bib{MR805809}{article}{
      author={Chanillo, Sagun},
      author={Wheeden, Richard~L.},
       title={Weighted {P}oincar\'e and {S}obolev inequalities and estimates
  for weighted {P}eano maximal functions},
        date={1985},
        ISSN={0002-9327},
     journal={Amer. J. Math.},
      volume={107},
      number={5},
       pages={1191\ndash 1226},
         url={http://dx.doi.org/10.2307/2374351},
      review={\MR{805809 (87f:42045)}},
}

\bib{MR2345338}{article}{
      author={Cohn, William~S.},
      author={Lu, Guozhen},
      author={Wang, Peiyong},
       title={Sub-elliptic global high order {P}oincar\'e inequalities in
  stratified {L}ie groups and applications},
        date={2007},
        ISSN={0022-1236},
     journal={J. Funct. Anal.},
      volume={249},
      number={2},
       pages={393\ndash 424},
         url={http://dx.doi.org/10.1016/j.jfa.2007.02.007},
      review={\MR{2345338 (2009j:43008)}},
}

\bib{MR0367477}{article}{
      author={Folland, G.~B.},
      author={Stein, Elias~M.},
       title={Estimates for the {$\bar \partial _{b}$} complex and analysis on
  the {H}eisenberg group},
        date={1974},
        ISSN={0010-3640},
     journal={Comm. Pure Appl. Math.},
      volume={27},
       pages={429\ndash 522},
      review={\MR{MR0367477 (51 \#3719)}},
}

\bib{MR1314734}{article}{
      author={Franchi, Bruno},
      author={Gallot, Sylvestre},
      author={Wheeden, Richard~L.},
       title={Sobolev and isoperimetric inequalities for degenerate metrics},
        date={1994},
        ISSN={0025-5831},
     journal={Math. Ann.},
      volume={300},
      number={4},
       pages={557\ndash 571},
         url={http://dx.doi.org/10.1007/BF01450501},
      review={\MR{1314734 (96a:46066)}},
}

\bib{MR1343563}{article}{
      author={Franchi, Bruno},
      author={Lu, Guozhen},
      author={Wheeden, Richard~L.},
       title={Representation formulas and weighted {P}oincar\'e inequalities
  for {H}\"ormander vector fields},
        date={1995},
        ISSN={0373-0956},
     journal={Ann. Inst. Fourier (Grenoble)},
      volume={45},
      number={2},
       pages={577\ndash 604},
         url={http://www.numdam.org/item?id=AIF_1995__45_2_577_0},
      review={\MR{1343563 (96i:46037)}},
}

\bib{MR1354890}{article}{
      author={Franchi, Bruno},
      author={Lu, Guozhen},
      author={Wheeden, Richard~L.},
       title={Weighted {P}oincar\'e inequalities for {H}\"ormander vector
  fields and local regularity for a class of degenerate elliptic equations},
        date={1995},
        ISSN={0926-2601},
     journal={Potential Anal.},
      volume={4},
      number={4},
       pages={361\ndash 375},
         url={http://dx.doi.org/10.1007/BF01053453},
        note={Potential theory and degenerate partial differential operators
  (Parma)},
      review={\MR{1354890 (97e:35018)}},
}

\bib{MR1383947}{article}{
      author={Franchi, Bruno},
      author={Lu, Guozhen},
      author={Wheeden, Richard~L.},
       title={A relationship between {P}oincar\'e-type inequalities and
  representation formulas in spaces of homogeneous type},
        date={1996},
        ISSN={1073-7928},
     journal={Internat. Math. Res. Notices},
      number={1},
       pages={1\ndash 14},
         url={http://dx.doi.org/10.1155/S1073792896000013},
      review={\MR{1383947 (97k:26012)}},
}

\bib{MR1744780}{article}{
      author={Franchi, Bruno},
      author={P{\'e}rez, Carlos},
      author={Wheeden, Richard~L.},
       title={Sharp geometric {P}oincar\'e inequalities for vector fields and
  non-doubling measures},
        date={2000},
        ISSN={0024-6115},
     journal={Proc. London Math. Soc. (3)},
      volume={80},
      number={3},
       pages={665\ndash 689},
         url={http://dx.doi.org/10.1112/S0024611500012375},
      review={\MR{1744780 (2001i:46048)}},
}

\bib{MR1455453}{article}{
   author={Franchi, Bruno},
   author={Wheeden, Richard L.},
   title={Compensation couples and isoperimetric estimates for vector
   fields},
   journal={Colloq. Math.},
   volume={74},
   date={1997},
   number={1},
   pages={9--27},
   issn={0010-1354},
   review={\MR{1455453 (98g:46042)}},
}


\bib{MR1421823}{incollection}{
      author={Gromov, Mikhael},
       title={Carnot-{C}arath\'eodory spaces seen from within},
        date={1996},
   booktitle={Sub-{R}iemannian geometry},
      series={Progr. Math.},
      volume={144},
   publisher={Birkh\"auser},
     address={Basel},
       pages={79\ndash 323},
      review={\MR{MR1421823 (2000f:53034)}},
}

\bib{MR1886617}{incollection}{
      author={Haj{\l}asz, Piotr},
       title={Sobolev inequalities, truncation method, and {J}ohn domains},
        date={2001},
   booktitle={Papers on analysis},
      series={Rep. Univ. Jyv\"askyl\"a Dep. Math. Stat.},
      volume={83},
   publisher={Univ. Jyv\"askyl\"a, Jyv\"askyl\"a},
       pages={109\ndash 126},
      review={\MR{1886617 (2003a:46052)}},
}

\bib{MR1800917}{book}{
      author={Heinonen, Juha},
       title={Lectures on analysis on metric spaces},
      series={Universitext},
   publisher={Springer-Verlag},
     address={New York},
        date={2001},
        ISBN={0-387-95104-0},
      review={\MR{MR1800917 (2002c:30028)}},
}

\bib{MR850547}{article}{
      author={Jerison, David},
       title={The {P}oincar\'e inequality for vector fields satisfying
  {H}\"ormander's condition},
        date={1986},
        ISSN={0012-7094},
     journal={Duke Math. J.},
      volume={53},
      number={2},
       pages={503\ndash 523},
      review={\MR{MR850547 (87i:35027)}},
}

\bib{MR1187073}{incollection}{
      author={Long, Rui~Lin},
      author={Nie, Fu~Sheng},
       title={Weighted {S}obolev inequality and eigenvalue estimates of
  {S}chr\"odinger operators},
        date={1991},
   booktitle={Harmonic analysis ({T}ianjin, 1988)},
      series={Lecture Notes in Math.},
      volume={1494},
   publisher={Springer, Berlin},
       pages={131\ndash 141},
         url={http://dx.doi.org/10.1007/BFb0087765},
      review={\MR{1187073 (94c:46066)}},
}

\bib{MR1202416}{article}{
      author={Lu, Guozhen},
       title={Weighted {P}oincar\'e and {S}obolev inequalities for vector
  fields satisfying {H}\"ormander's condition and applications},
        date={1992},
        ISSN={0213-2230},
     journal={Rev. Mat. Iberoamericana},
      volume={8},
      number={3},
       pages={367\ndash 439},
         url={http://dx.doi.org/10.4171/RMI/129},
      review={\MR{1202416 (94c:35061)}},
}

\bib{MR1286482}{article}{
      author={Lu, Guozhen},
       title={The sharp {P}oincar\'e inequality for free vector fields: an
  endpoint result},
        date={1994},
        ISSN={0213-2230},
     journal={Rev. Mat. Iberoamericana},
      volume={10},
      number={2},
       pages={453\ndash 466},
         url={http://dx.doi.org/10.4171/RMI/158},
      review={\MR{1286482 (96g:26023)}},
}

\bib{MR1642822}{article}{
      author={Lu, Guozhen},
      author={Wheeden, Richard~L.},
       title={An optimal representation formula for {C}arnot-{C}arath\'eodory
  vector fields},
        date={1998},
        ISSN={0024-6093},
     journal={Bull. London Math. Soc.},
      volume={30},
      number={6},
       pages={578\ndash 584},
         url={http://dx.doi.org/10.1112/S0024609398004895},
      review={\MR{1642822 (2000a:31005)}},
}

\bib{MR1792599}{article}{
      author={Lu, Guozhen},
      author={Wheeden, Richard~L.},
       title={High order representation formulas and embedding theorems on
  stratified groups and generalizations},
        date={2000},
        ISSN={0039-3223},
     journal={Studia Math.},
      volume={142},
      number={2},
       pages={101\ndash 133},
      review={\MR{1792599 (2001k:46055)}},
}

\bib{MR0340523}{article}{
      author={Muckenhoupt, Benjamin},
      author={Wheeden, Richard~L.},
       title={Weighted norm inequalities for fractional integrals},
        date={1974},
        ISSN={0002-9947},
     journal={Trans. Amer. Math. Soc.},
      volume={192},
       pages={261\ndash 274},
      review={\MR{0340523 (49 \#5275)}},
}

\bib{MR1882665}{article}{
      author={Nagel, Alexander},
      author={Stein, Elias~M.},
       title={Differentiable control metrics and scaled bump functions},
        date={2001},
        ISSN={0022-040X},
     journal={J. Differential Geom.},
      volume={57},
      number={3},
       pages={465\ndash 492},
         url={http://projecteuclid.org/getRecord?id=euclid.jdg/1090348130},
      review={\MR{MR1882665 (2003i:58003)}},
}

\bib{MR793239}{article}{
      author={Nagel, Alexander},
      author={Stein, Elias~M.},
      author={Wainger, Stephen},
       title={Balls and metrics defined by vector fields. {I}. {B}asic
  properties},
        date={1985},
        ISSN={0001-5962},
     journal={Acta Math.},
      volume={155},
      number={1-2},
       pages={103\ndash 147},
      review={\MR{MR793239 (86k:46049)}},
}

\bib{MR1818113}{article}{
      author={P{\'e}rez, Carlos},
      author={Wheeden, Richard~L.},
       title={Uncertainty principle estimates for vector fields},
        date={2001},
        ISSN={0022-1236},
     journal={J. Funct. Anal.},
      volume={181},
      number={1},
       pages={146\ndash 188},
         url={http://dx.doi.org/10.1006/jfan.2000.3711},
      review={\MR{1818113 (2002h:42035)}},
}

\bib{MR1962949}{article}{
      author={P{\'e}rez, Carlos},
      author={Wheeden, Richard~L.},
       title={Potential operators, maximal functions, and generalizations of
  {$A_\infty$}},
        date={2003},
        ISSN={0926-2601},
     journal={Potential Anal.},
      volume={19},
      number={1},
       pages={1\ndash 33},
         url={http://dx.doi.org/10.1023/A:1022449810008},
      review={\MR{1962949 (2004a:42025)}},
}

\bib{MR0436223}{article}{
      author={Rothschild, Linda~Preiss},
      author={Stein, Elias~M.},
       title={Hypoelliptic differential operators and nilpotent groups},
        date={1976},
        ISSN={0001-5962},
     journal={Acta Math.},
      volume={137},
      number={3-4},
       pages={247\ndash 320},
      review={\MR{MR0436223 (55 \#9171)}},
}

\bib{MR762360}{article}{
      author={S{\'a}nchez-Calle, Antonio},
       title={Fundamental solutions and geometry of the sum of squares of
  vector fields},
        date={1984},
        ISSN={0020-9910},
     journal={Invent. Math.},
      volume={78},
      number={1},
       pages={143\ndash 160},
         url={http://dx.doi.org/10.1007/BF01388721},
      review={\MR{762360 (86e:58078)}},
}

\bib{MR1175693}{article}{
      author={Sawyer, Eric},
      author={Wheeden, Richard~L.},
       title={Weighted inequalities for fractional integrals on {E}uclidean and
  homogeneous spaces},
        date={1992},
        ISSN={0002-9327},
     journal={Amer. J. Math.},
      volume={114},
      number={4},
       pages={813\ndash 874},
         url={http://dx.doi.org/10.2307/2374799},
      review={\MR{1175693 (94i:42024)}},
}

\bib{MR1070036}{article}{
      author={Varopoulos, N.~Th.},
       title={Small time {G}aussian estimates of heat diffusion kernels. {II}.
  {T}he theory of large deviations},
        date={1990},
        ISSN={0022-1236},
     journal={J. Funct. Anal.},
      volume={93},
      number={1},
       pages={1\ndash 33},
         url={http://dx.doi.org/10.1016/0022-1236(90)90136-9},
      review={\MR{MR1070036 (91j:58156)}},
}

\bib{PL}{article}{
      author={Yung, Po-Lam},
       title={Sobolev inequalities for $(0,q)$ forms on {C}{R} manifolds of
  finite type},
     journal={Math. Res. Letts., 17 (2010), no. 1, 177-196},
}

\end{biblist}
\end{bibdiv}

\end{document}